\documentclass[aos]{imsart}

\RequirePackage{amsthm,amsmath,amsfonts,amssymb}
\RequirePackage[numbers]{natbib}
\RequirePackage[colorlinks,citecolor=blue,urlcolor=blue]{hyperref}
\RequirePackage{graphicx}

\RequirePackage{algorithm}
\RequirePackage{algpseudocode}

\RequirePackage{centernot}

\startlocaldefs
\theoremstyle{plain}

\newtheorem{theorem}{Theorem}[section]
\newtheorem{proposition}{Proposition}[section]
\newtheorem{lemma}[theorem]{Lemma}
\theoremstyle{remark}
\newtheorem{definition}[theorem]{Definition}
\newtheorem{example}{Example}

\endlocaldefs

\begin{document}

\begin{frontmatter}
\title{On the Independence Polynomial and Threshold of an Antiregular $k$-Hypergraph}
\runtitle{Independence Polynomials and Thresholds of Antiregular $k$-Hypergraphs}

\begin{aug}
\author[A,B]{\fnms{Erchuan}~\snm{Zhang}\ead[label=e2]{erchuan.zhang@ecu.edu.au}\ead[label=e3]{erchuan.zhang@uwa.edu.au}\orcid{0000-0002-4005-5431}}
\address[A]{School of Science,
Edith Cowan University\printead[presep={, \ }]{e2}}

\address[B]{Department of Mathematics and Statistics,
The University of Western Australia\printead[presep={, \ }]{e3}}
\end{aug}

\begin{abstract}
Given an integer $k\geq 3$ and an initial $k-1$ isolated vertices, an {\em antiregular $k$-hypergraph} is constructed by alternatively adding an isolated vertex (connected to no other vertices) or a dominating vertex (connected to every other $k-1$ vertices). Let $a_i$ be the number of independent sets of cardinality $i$ in a hypergraph $H$, then the {\em independence polynomial} of $H$ is defined as $I(H;x)=\sum_{i=0}^m a_i x^i$, where $m$ is the size of a maximum independent set. The main purpose of the present paper is to generalise some results of independence polynomials of antiregular graphs to the case of antiregular $k$-hypergraphs. In particular, we derive (semi-)closed formulas for the independence polynomials of antiregular $k$-hypergraphs and prove their log-concavity. Furthermore, we show that antiregular $k$-hypergraphs are {\em $T2$-threshold}, which means there exist a labeling $c$ of the vertex set and a threshold $\tau$ such that for any vertex subset $S$ of cardinality $k$, $\sum_{i\in S}c(i)>\tau$ if and only if $S$ is a hyperedge.
\end{abstract}



\end{frontmatter}

\section{Introduction}

Let $G=(V,E)$ be a simple (i.e., finite, undirected, loopless, no multiple edges) graph with the vertex set $V$ and edge set $E$, then $G$ is called {\em antiregular} \cite{merris2003antiregular}, {\em quasiperfect} \cite{behzad1967no}, {\em maximally nonregular} \cite{zykov1990fundamentals} if its vertex degrees take on $\vert V\vert-1$ different values, i.e., two vertices share the same degree, where $\vert V\vert$ is the cardinality of $V$.

Let $G_1,G_2$ be two simple graphs, the {\em disjoint union} of them, denoted by $G_1\cup G_2$, is the graph $G$ with the vertex set $V(G)=V(G_1)\cup V(G_2)$ and the edge set $E(G)=E(G_1)\cup E(G_2)$. For two given disjoint graphs $G_1$ and $G_2$, their {\em Zykov sum}, denoted by $G_1+G_2$, is the graph $G$ with the vertex set $V(G)=V(G_1)\cup V(G_2)$ and the edge set $E(G)=E(G_1)\cup E(G_2)\cup\{uv\vert u\in V(G_1),v\in V(G_2)\}$. Then the antiregular graphs can also be defined by the following recurrence relationship \cite{merris2003antiregular}
\begin{equation}\label{rel1}
\begin{aligned}
&A_1=K_1, A_{n+1}=K_1+\bar{A}_n, n\geq 1,~\hbox{or}\\
&A_1=K_1, A_2=K_2, A_{n+2}=K_1+(K_1\cup A_n), n\geq 1,
\end{aligned}
\end{equation}
where $K_1$ denotes a graph with an isolated vertex only, $K_2$ denotes the complete graph on two vertices, $\bar{A}_n$ is the complement of $A_n$.

By the relation \eqref{rel1}, an antiregular graph can be constructed by adding an isolated vertex or a dominating vertex alternatively. A vertex is called {\em isolated} if it is connected to no other vertex, and {\em dominating} if it is connected to every other vertex. Then, an antiregular graph can be represented by a binary string, called {\em binary building string} \cite{levit2012independence},
\begin{align}
{\bf 00}101010\ldots ~~\hbox{or}~~ {\bf 01}010101\ldots,
\end{align}
where $0$ stands for adding an isolated vertex and $1$ stands for adding a dominating vertex. Recall that a graph $G$ is called {\em threshold} if it can be constructed from $K_1$ by iterating the operations of adding an isolated vertex and adding a dominating vertex in any order \cite{chvatal1973set}. Therefore, every antiregular graph is threshold.

Given a simple graph $G=(V,E)$, a vertex subset $W\subseteq V$ is called an {\em independent set} if it does not include two adjacent vertices. The {\em independence polynomial} $I(G;x)$ of a graph $G$ is defined as \cite{gutman1983generalizations}
\begin{align}
I(G;x)=\sum_{i=0}^\alpha s_ix^i,
\end{align}
where $s_i$ is the number of independent sets of cardinality $i$ in the graph $G$, $0\leq i\leq \alpha$, $\alpha$ is the size of a maximum independent set.

Levit and Mandrescu in \cite{levit2012independence} derived closed formulas for the independence polynomial $I(G;x)$ of an antiregular graph $G$ and showed that $I(G;x)$ is log-concave with at most two real roots. Moreover, they found that antiregular graph is uniquely determined by its independence polynomial within the family of threshold graphs, which follows from the fact that independence polynomials of threshold graphs are unique \cite{hoede1994clique}.

The motivation of the present paper is to generalise some results of independence polynomials of antiregular graphs in \cite{levit2012independence} to the so-called {\em antiregular $k$-hypergraphs}. A {\em hypergraph} $H=(V,E)$ is an order pair of a set of vertices $V$ and a set of hyperedges $E$, where each hyperedge $e\in E$ is a non-empty subset of the vertex set $V$. Further, if $\vert e\vert=k$ for any $e\in E$, then $H$ is called {\em k-uniform hypergraph}, {\em k-hypergraph} or {\em k-families}. Note that $2$-hypergraphs are the usual graphs. Throughout this paper, we set $k\geq 3$ if no specific declaration.

The {\em disjoint union} of two hypergraphs $H_1,H_2$ is the hypergraph $H=H_1\cup H_2$ having the disjoint union of $V(H_1)$, $V(H_2)$ as a vertex set, and the disjoint union of $E(H_1)$, $E(H_2)$ as a hyperedge set. Let $H_1$, $H_2$ be disjoint hypergraphs and $\vert V(H_2)\vert\geq k-1$ for some $k\geq 2$, their {\em generalized Zykov k-sum} is the hypergraph $H=H_1\oplus_k H_2$ with $V(H_1)\cup V(H_2)$ as a vertex set and $E(H_1)\cup E(H_2)\cup \{v w_1 w_2\cdots w_{k-1}\vert v\in V(H_1),w_1,\cdots,w_{k-1}\in V(H_2)\}$ as a hyperedge set. Note that the generalised Zykov $2$-sum is the usual Zykov sum. For $k>2$, the generalised Zykov $k$-sum is not commutative in general, i.e., $H_1\oplus_k H_2\neq H_2\oplus_k H_1$ for general hypergraphs $H_1,H_2$.

Now we give a formal definition of an antiregular $k$-hypergraph as follows.

\begin{definition}\label{antidef}
For $k\geq 3$, the antiregular $k$-hypergraphs can be defined by the following recurrences,
\begin{equation}
\begin{aligned}
&A_i=i K_1, 1\leq i\leq k-1,\\
&A_{n+2}=K_1\oplus_k \bar{A}_{n+1}=K_1\oplus_k (K_1\cup A_n), \hbox{or}\\
&\bar{A}_{n+2}=K_1\cup A_{n+1}=K_1\cup (K_1\oplus_k \bar{A}_n), n\geq k-1,
\end{aligned}
\end{equation}
where $K_1$ is the hypergraph with one isolated vertex only, $i K_1$ denotes the disjoint union of $i$ copies of the hypergraph $K_1$, $1\leq i\leq k-1$, $\bar{A}_n$ is the complement of $A_n$. We call $A_n$ the {\em connected} antiregular $k$-hypergraph with $n$ vertices, $\bar{A}_n$ the {\em disconnected} one. 
\end{definition}

Similar to the graph case, we use $0$ to denote adding an isolated vertex (connected to no other vertices) and $1$ to denote adding a dominating vertex (connected to every other $k-1$ vertices). Then the antiregular $k$-hypergraphs can be represented by the following binary building strings,
\begin{equation}
\begin{aligned}
&(k=3)~~{\bf 000}101010\ldots ~~\hbox{or}~~ {\bf 001}010101\ldots,\\
&(k> 3)~~0[k]101010\ldots ~~\hbox{or}~~ 0[k-1]1010101\ldots,
\end{aligned}
\end{equation}
where $0[k]$ means a string with $k$ zeros. If the last bit of the binary building string of an antiregular $k$-hypergraph is $1$, then the hypergraph is connected. Otherwise, it is disconnected. More generally, we call a $k$-hypergraph {\em $\{0,1\}$-constructable} if it can be constructed by adding isolated and dominating vertices in some order.

The {\em vertex-degree} of a vertex in a hypergraph is the number of hyperedges that contain this vertex. By Definition \ref{antidef}, the first $k$ added vertices have the same vertex-degree. Moreover, from the $(k+1)$-th or $(k+2)$-th added isolated vertex to the last one, their vertex-degrees are decreasing and from the first added dominating vertex to the last one, their vertex-degrees are increasing. Therefore, we have the following proposition.

\begin{proposition}
An antiregular $k$-hypergraph only has $k$ vertices with the same vertex-degree.
\end{proposition}

\begin{proof}
Let $b^n=0[k-1]1010\ldots$ be the binary building string of an antiregular $k$-hypergraph with $n$ vertices $(n\geq k)$ (It is similar to show for $b^n=0[k]1010\ldots$). We will show that the vertex-degree $d$ satisfies
\begin{align}\label{vertexpropeq1}
d(1)=\ldots=d(k)>d(i_1)>d(i_2)
\end{align}
for any isolated vertices $k<i_1<i_2$, and 
\begin{align}\label{vertexpropeq2}
d(j_1)<d(j_2)
\end{align}
for any dominating vertices $j_1<j_2$ by induction on $n$.

If $n=k$, then $d(1)=\ldots=d(k)=1$. If $n=k+1$, then $d(1)=\ldots=d(k)>d(k+1)=0$. If $n=k+2$, then $d(i)=1+\binom{k}{k-2}$, where $1\leq i\leq k$, and $d(k+1)=\binom{k}{k-2}$, $d(k+2)=\binom{k+1}{k-1}$. Thus, $d(k-1)>d(k+1)$ and $d(k)<d(k+2)$.

Suppose \eqref{vertexpropeq1} and \eqref{vertexpropeq2} hold for $n\leq m$. Let $d^\prime,d$ be the vertex-degree sequences of $b^m,b^{m+1}$, respectively. If $m+1$ is an isolated vertex, then $d(m+1)=0$, $d(i)=d^\prime(i)$, where $1\leq i\leq m$. Thus, $d(m-1)>d(m+1)$. If $m+1$ is a dominating vertex, then $d(m+1)=\binom{m}{k-1}$ and $d(i)=d^\prime(i)+\binom{m-1}{k-2}$, where $1\leq i\leq m$. Thus, $d(m-1)=\binom{m-2}{k-1}+\binom{m-1}{k-2}<d(m+1)$.
\end{proof}

In simple graphs, independent sets contain no two adjacent vertices, i.e., contain no edges. Motivated by this idea, a vertex subset $W\subseteq V$ is called an {\em independent set} in a hypergraph $H=(V,E)$ if $W$ does not include any hyperedges, i.e., $\forall e\in E$, $e\nsubseteq W$. Given a hypergraph $H=(V,E)$, the {\em independence polynomial} $I(H;x)$ of $H$ is defined as \cite{trinks2016survey}
\begin{align}
I(H;x)=\sum_{W\subseteq V~\hbox{is independent in}~H} x^{\vert W\vert}.
\end{align}
For $k$-hypergraphs $H$, any vertex subset of size less than $k$ cannot contain a hyperedge, thus, the coefficient of $x^i$ in $I(H;x)$ is $\binom{\vert V(H)\vert}{i}$ for $0\leq i\leq k-1$, and the coefficient of $x^k$ is $\binom{\vert V(H)\vert}{k}-\vert E(H)\vert$.


The main contributions of the present paper are listed as follows.
\begin{itemize}
\item We derive the recurrence relations (Theorem \ref{recthm}) and (semi-)closed forms (Theorems \ref{k3ipthm}, \ref{kipthm}) of independence polynomials of antiregular $k$-hypergraphs.
\item We show that independence polynomials of antiregular $k$-hypergraphs are log-concave (Theorems \ref{logc3thm}, \ref{logckthm}).
\item Every $\{0,1\}$-constructable (including antiregular) $k$-hypergraph is shown to be $T2$-threshold (Theorem \ref{thresholdthm}), which means there exist a labeling $c$ of the vertex set and a threshold $\tau$ such that for any vertex subset $S$ of cardinality $k$, $\sum_{i\in S}c(i)>\tau$ if and only if $S$ is a hyperedge.
\end{itemize}

The remainder of this paper is organised as follows. In section \ref{sec2}, with the help of a general recurrence relationship of independence polynomial of hypergraphs derived by Trinks \cite{trinks2016survey}, we present the recurrence relations of independence polynomials of antiregular $k$-hypergraphs. Furthermore, (semi-)closed formulas for independence polynomials are obtained. In section \ref{sec3}, by the induction on the number of vertices in an antiregular $k$-hypergraph, we show that independence polynomials of antiregular $k$-hypergraph are log-concave. In section \ref{sec4}, we develop an algorithm and use it to prove that $\{0,1\}$-constructable (including antiregular) $k$-hypergraphs are $T2$-threshold. A conclusion is given in the following section.

\section{Closed and semi-closed forms of independence polynomials}\label{sec2}

Given a hypergraph $H=(V,E)$ and a vertex subset $W\subseteq V$, we define the following hypergraph operations:
\begin{itemize}
\item {\em deletion} of the vertices $w\in W$, denoted by $\ominus W$, means the vertices $w$ and their incident hyperedges are removed.
\item {\em hiding} of the vertices $w\in W$, denoted by $\sim W$, means the vertices $w$ are removed in the vertex set $V$ and in their incident hyperedges.
\end{itemize}

The following recurrence relations of independence polynomials of hypergraphs are established by Trinks \cite[Theorem 3]{trinks2016survey}.

\begin{lemma}\label{reclem}
Let $H=(V,E)$ be a hypergraph and $v\in V$, the independence polynomial $I(H;x)$ of $H$ satisfies
\begin{align}
I(H;x)=\begin{cases}
I(H_{\ominus v};x)+x I(H_{\sim v};x)~~~~&\hbox{if}~~\{v\}\notin E,\\
I(H_{\ominus v};x) &\hbox{otherwise.}
\end{cases}
\end{align}
\end{lemma}

Based on the definition of antiregular $k$-hypergraphs and Lemma \ref{reclem}, we can now derive the recurrence relations of independence polynomials for antiregular $k$-hypergraphs specifically.

\begin{theorem}\label{recthm}
Let $A_n$ be an antiregular $k$-hypergraph with $n$ vertices and $\bar{A}_n$ the complement of $A_n$, $n\geq 1$, its independence polynomial satisfies the following recurrence relationship
\begin{align}
&\begin{cases}
I(A_i;x)=(1+x)^i,~~ 1\leq i\leq k-1,\\
I(\bar{A}_{n+1};x)=(1+x)I(A_n;x),\\
 I(A_{n+2};x)=I(\bar{A}_{n+1};x)+\sum_{i=1}^{k-1}\binom{n+1}{i-1}x^i, ~~n\geq k-1.
\end{cases}
\hbox{or}\\
&\begin{cases}
I(A_i;x)=(1+x)^i, ~~1\leq i\leq k-1,\\
I(A_{n+1};x)=I(\bar{A}_{n};x)+\sum_{i=1}^{k-1}\binom{n}{i-1}x^i,\\
I(\bar{A}_{n+2};x)=(1+x)I(A_{n+1};x), ~~n\geq k-1.
\end{cases}
\end{align}
\end{theorem}

\begin{proof}
For $1\leq i\leq k-1$, $I(A_i;x)=I(iK_1;x)=(I(K_1;x))^i=(1+x)^i$. We now only prove the case where $A_{n+2}=K_1\oplus_k\bar{A}_{n+1}=K_1\oplus_k(K_1\cup A_n)$ as the other case is similar to show.

Since $\bar{A}_{n+1}=K_1\cup A_n$, we have
\begin{align*}
\bar{A}_{n+1 \ominus K_1}=\bar{A}_{n+1 \sim K_1}=A_n.
\end{align*}
By Lemma \ref{reclem}, we get
\begin{align*}
I(\bar{A}_{n+1};x)=I(\bar{A}_{n+1 \ominus K_1};x)+x I(\bar{A}_{n+1 \sim K_1};x)=(1+x)I(A_n;x).
\end{align*}

Since $A_{n+2}=K_1\oplus_k\bar{A}_{n+1}$, we have $A_{n+2 \ominus K_1}=\bar{A}_{n+1}$ and
\begin{align*}
I(A_{n+2\sim K_1};x)=\sum_{i=0}^{k-2}\binom{n+1}{i}x^i,
\end{align*}
which is because any vertex subset of size $k-1$ of $A_{n+2 \sim K_1}$ is a hyperedge. By Lemma \ref{reclem}, we obtain
\begin{align*}
I(A_{n+2};x)=I(A_{n+2 \ominus K_1};x)+x I(A_{n+2 \sim K_1};x)=I(\bar{A}_{n+1};x)+\sum_{i=1}^{k-1}\binom{n+1}{i-1}x^i.
\end{align*}
\end{proof}

\begin{theorem}\label{k3ipthm}
The independence polynomials of antiregular $3$-hypergraphs $A_n$ and $\bar{A}_n$ are given by the following closed forms,
\begin{align}
&\begin{cases}
I(A_{2n-1};x)=3(1+x)^n+(1+x)^{n-1}-2nx-3,\\
I(\bar{A}_{2n};x)=3(1+x)^{n+1}+(1+x)^n-(1+x)(2nx+3),
\end{cases}
n\geq 1, ~~\hbox{or}\label{k3ipeq1}\\
&\begin{cases}
I(\bar{A}_{2n-1};x)=(1+x)^{n+1}+3(1+x)^n-(1+x)((2n-1)x+3),\\
I(A_{2n};x)=(1+x)^{n+1}+3(1+x)^n-(2n+1)x-3,
\end{cases}
n\geq 1.
\end{align}
\end{theorem}

\begin{proof}
We only prove the formula \eqref{k3ipeq1} by induction on $n$ and the other one is similar to show.

If $n=1$, then
\begin{align*}
&I(A_1;x)=I(K_1;x)=1+x=3(1+x)+1-2x-3,\\
&I(\bar{A}_2;x)=I(2K_1;x)=(1+x)^2=3(1+x)^2+1+x-(1+x)(2x+3).
\end{align*}

Suppose the formula \eqref{k3ipeq1} is true for $n=1,\cdots,m$.

If $m+1=2s+1$, then, by Theorem \ref{recthm} and induction hypothesis, we have
\begin{align*}
I(A_{2s+1};x)&=I(\bar{A}_{2s};x)+x(2sx+1)\\
&=3(1+x)^{s+1}+(1+x)^s-(1+x)(2sx+3)+x(2sx+1)\\
&=3(1+x)^{s+1}+(1+x)^s-2(s+1)x-3.
\end{align*}

If $m+1=2s+2$, using Theorem \ref{recthm} and induction hypothesis again, we find
\begin{align*}
I(\bar{A}_{2s+2};x)&=(1+x)I(A_{2s+1};x)\\
&=(1+x)\left(3(1+x)^{s+1}+(1+x)^s-2(s+1)x-3\right)\\
&=3(1+x)^{s+2}+(1+x)^{s+1}-(1+x)(2(s+1)x+3).
\end{align*}
Therefore, this theorem is true by induction.
\end{proof}

From the theorem above, we can observe that once the formulas of independence polynomials are given, it is not hard to check their correctness by the induction method. In what follows, we will derive the formulas of independence polynomials for general antiregular $k$-hypergraphs from Theorem \ref{recthm} directly.

\begin{theorem}\label{kipthm}
The independence polynomials of antiregular $k$-hypergraphs $A_n$ and $\bar{A}_n$ are given by,
for $n\geq \frac{k+1}{2}$,
\begin{align}
&I(A_{2n-1};x)=\begin{cases}
(1+x)^{n-1-\frac{k-1}{2}}\left[(1+x)^{k-1}+\sum_{i=0}^{k-1}\alpha_i^{k-1} x^i\right]-\sum_{i=0}^{k-1}\alpha_i^{2n-2}x^i\\
+\sum_{i=1}^{k-1}\binom{2n-2}{i-1}x^i,~~k~\hbox{is odd,}\\
(1+x)^{n-1-\frac{k}{2}}\left[(1+x)^{k}+\sum_{i=0}^{k-1}\alpha_i^{k} x^i\right]-\sum_{i=0}^{k-1}\alpha_i^{2n-2}x^i\\
+\sum_{i=1}^{k-1}\binom{2n-2}{i-1}x^i,~~k~\hbox{is even,}
\end{cases}\label{kipeq1}\\
&I(\bar{A}_{2n};x)=\begin{cases}
(1+x)^{n-\frac{k-1}{2}}\left[(1+x)^{k-1}+\sum_{i=0}^{k-1}\alpha_i^{k-1} x^i\right]-\sum_{i=0}^{k-1}\alpha_i^{2n}x^i,~~k~\hbox{is odd,}\\
(1+x)^{n-\frac{k}{2}}\left[(1+x)^{k}+\sum_{i=0}^{k-1}\alpha_i^{k} x^i\right]-\sum_{i=0}^{k-1}\alpha_i^{2n}x^i,~~k~\hbox{is even,}
\end{cases}\label{kipeq2}
\end{align}
where $\{\alpha_i^{2n}\}$ satisfy
\begin{align}\label{kipeq3}
\begin{cases}
\alpha_0^{2n-2}-\alpha_0^{2n}=0,\\
\alpha_i^{2n-2}+\alpha_{i-1}^{2n-2}-\alpha_i^{2n}=\binom{2n-1}{i-1},~~1\leq i\leq k-1,\\
\alpha_{k-1}^{2n-2}=\binom{2n-2}{k-2}.
\end{cases}
\end{align}
Or, 
\begin{align}
&I(\bar{A}_{2n-1};x)=\begin{cases}
(1+x)^{n-1-\frac{k-1}{2}}\left[(1+x)^{k}+\sum_{i=0}^{k-1}\beta_i^{k} x^i\right]-\sum_{i=0}^{k-1}\beta_i^{2n-1}x^i,~~k~\hbox{is odd,}\\
(1+x)^{n-\frac{k}{2}}\left[(1+x)^{k-1}+\sum_{i=0}^{k-1}\beta_i^{k-1} x^i\right]-\sum_{i=0}^{k-1}\beta_i^{2n-1}x^i,~~k~\hbox{is even,}
\end{cases}\label{kipeq4}\\
&I(A_{2n};x)=\begin{cases}
(1+x)^{n-1-\frac{k-1}{2}}\left[(1+x)^{k}+\sum_{i=0}^{k-1}\beta_i^{k} x^i\right]-\sum_{i=0}^{k-1}\beta_i^{2n-1}x^i\\
+\sum_{i=1}^{k-1}\binom{2n-1}{i-1}x^i,~~k~\hbox{is odd,}\\
(1+x)^{n-\frac{k}{2}}\left[(1+x)^{k-1}+\sum_{i=0}^{k-1}\beta_i^{k-1} x^i\right]-\sum_{i=0}^{k-1}\beta_i^{2n-1}x^i\\
+\sum_{i=1}^{k-1}\binom{2n-1}{i-1}x^i,~~k~\hbox{is even,}
\end{cases}\label{kipeq5}
\end{align}
where $\{\beta_i^{2n}\}$ satisfy
\begin{align}\label{kipeq6}
\begin{cases}
\beta_0^{2n-1}-\beta_0^{2n+1}=0,\\
\beta_i^{2n-1}+\beta_{i-1}^{2n-1}-\beta_i^{2n+1}=\binom{2n}{i-1},~~1\leq i\leq k-1,\\
\beta_{k-1}^{2n-1}=\binom{2n-1}{k-2}.
\end{cases}
\end{align}
\end{theorem}

\begin{proof}
We only show the formula \eqref{kipeq2}, and \eqref{kipeq1} follows from \eqref{kipeq2} and $I(A_{2n-1};x)=I(\bar{A}_{2n-2};x)+\sum_{i=1}^{k-1}\binom{2n-2}{i-1}x^i$ directly. The proof for \eqref{kipeq4} and \eqref{kipeq5} is similar.

By Theorem \ref{recthm}, we have
\begin{equation}\label{kipeq7}
\begin{aligned}
I(\bar{A}_{2n};x)&=(1+x)I(A_{2n-1};x)\\
&=(1+x)\left[I(\bar{A}_{2n-2};x)+\sum_{i=1}^{k-1}\binom{2n-2}{i-1}x^i\right]\\
&=(1+x)I(\bar{A}_{2n-2};x)+\sum_{i=1}^{k-1}\binom{2n-1}{i-1}x^i+\binom{2n-2}{k-2}x^k,
\end{aligned}
\end{equation}
where we have used the formula $\binom{2n-2}{i-1}+\binom{2n-2}{i-2}=\binom{2n-1}{i-1}$.

Suppose
\begin{align}
I(\bar{A}_{2n};x)+\sum_{i=0}^{k-1}\alpha_i^{2n}x^i=(1+x)\left[I(\bar{A}_{2n-2};x)+\sum_{i=0}^{k-1}\alpha_i^{2n-2}x^i\right],
\end{align}
which implies
\begin{equation}\label{kipeq8}
\begin{aligned}
I(\bar{A}_{2n};x)&=(1+x)I(\bar{A}_{2n-2};x)+\sum_{i=0}^{k-1}\left(\alpha_i^{2n-2}-\alpha_i^{2n}\right)x^i+\sum_{i=1}^k\alpha_{i-1}^{2n-2}x^i\\
&=(1+x)I(\bar{A}_{2n-2};x)+\alpha_0^{2n-2}-\alpha_0^{2n}+\sum_{i=1}^{k-1}\left(\alpha_i^{2n-2}+\alpha_{i-1}^{2n-2}-\alpha_i^{2n}\right)x^i+\alpha_{k-1}^{2n-2}x^k.
\end{aligned}
\end{equation}
By comparing \eqref{kipeq7} and \eqref{kipeq8}, we find $\{\alpha_i^{2n}\}$ satisfy the relation \eqref{kipeq3}.

If $k$ is odd,
\begin{align*}
I(\bar{A}_{2n};x)+\sum_{i=0}^{k-1}\alpha_i^{2n}x^i&=(1+x)\left[I(\bar{A}_{2n-2};x)+\sum_{i=0}^{k-1}\alpha_i^{2n-2}x^i\right]\\
&=(1+x)^{n-\frac{k-1}{2}}\left[I(\bar{A}_{k-1};x)+\sum_{i=0}^{k-1}\alpha_i^{k-1}x^i\right].\\
\implies
I(\bar{A}_{2n};x)=&(1+x)^{n-\frac{k-1}{2}}\left[(1+x)^{k-1}+\sum_{i=0}^{k-1}\alpha_i^{k-1}x^i\right]-\sum_{i=0}^{k-1}\alpha_i^{2n}x^i.
\end{align*}

If $k$ is even,
\begin{align*}
I(\bar{A}_{2n};x)+\sum_{i=0}^{k-1}\alpha_i^{2n}x^i&=(1+x)\left[I(\bar{A}_{2n-2};x)+\sum_{i=0}^{k-1}\alpha_i^{2n-2}x^i\right]\\
&=(1+x)^{n-\frac{k}{2}}\left[I(\bar{A}_{k};x)+\sum_{i=0}^{k-1}\alpha_i^{k}x^i\right].\\
\implies
I(\bar{A}_{2n};x)=&(1+x)^{n-\frac{k}{2}}\left[(1+x)^k+\sum_{i=0}^{k-1}\alpha_i^{k}x^i\right]-\sum_{i=0}^{k-1}\alpha_i^{2n}x^i.
\end{align*}
\end{proof}

Though it is not hard to solve the algebraic equations \eqref{kipeq3} and \eqref{kipeq6} for small value of $k$, it would be very nasty to write down closed forms of $\alpha_i^{2n}$ and $\beta_i^{2n+1}$ for general $k$. Now we verify Theorem \ref{kipthm} for $k=2$ and $k=3$.

For $k=2$, solving \eqref{kipeq3},
\begin{align*}
\begin{cases}
\alpha_0^{2n-2}-\alpha_0^{2n}=0\\
\alpha_1^{2n-2}+\alpha_0^{2n-2}-\alpha_1^{2n}=1\\
\alpha_1^{2n-2}=1
\end{cases}
\implies 
\begin{cases}
\alpha_0^{2n}=1,\\
\alpha_1^{2n}=1.
\end{cases}
\end{align*}
Then, by \eqref{kipeq2}, we have
\begin{align*}
I(\bar{A}_{2n};x)&=(1+x)^{n-1}\left[(1+x)^2+\alpha_0^2+\alpha_1^2x\right]-\alpha_0^{2n}-\alpha_1^{2n}x\\
&=(1+x)^{n+1}+(1+x)^{n}-x-1,
\end{align*}
which has also been derived in \cite[Theorem 2.6]{levit2012independence}.

For $k=3$, solving \eqref{kipeq3}
\begin{align*}
\begin{cases}
\alpha_0^{2n-2}-\alpha_0^{2n}=0\\
\alpha_1^{2n-2}+\alpha_0^{2n-2}-\alpha_1^{2n}=1\\
\alpha_2^{2n-2}+\alpha_1^{2n-2}-\alpha_2^{2n}=2n-1\\
\alpha_2^{2n-2}=2n-2
\end{cases}
\implies
\begin{cases}
\alpha_0^{2n}=3,\\
\alpha_1^{2n}=2n+3,\\
\alpha_2^{2n}=2n.
\end{cases}
\end{align*}
Then, by \eqref{kipeq2}, we have
\begin{align*}
I(\bar{A}_{2n};x)&=(1+x)^{n-1}\left[(1+x)^2+\alpha_0^2+\alpha_1^2x+\alpha_2^2x^2\right]-\alpha_0^{2n}-\alpha_1^{2n}x-\alpha_2^{2n}x^2\\
&=3(1+x)^{n+1}+(1+x)^{n}-(1+x)(2nx+3),
\end{align*}
which is consistent with Theorem \ref{k3ipthm}.


\section{Log-concavity of independence polynomials}\label{sec3}

A finite sequence of real numbers $(a_1,a_2,\cdots,a_n)$ is called {\em log-concave} if $a_i^2\geq a_{i-1}a_{i+1}$ for $i=2,3,\cdots,n-1$. A polynomial is called {\em log-concave} if the sequence of its coefficients is log-concave. The product of two log-concave polynomials is log-concave \cite{keilson1971some}. In this section, we will show that the independence polynomial of an antiregular $k$-hypergraph is log-concave.

Let $G=A_n$ ($\bar{A}_n$) be a connected (disconnected) antiregular $k$-hypergraph with $n$ vertices. Suppose its independence polynomial is given by
\begin{align}
I(G;x)=\sum_{i=0}^m a_i^n x^i,
\end{align}
where $m$ is the size of a maximum independent set. Since any vertex subset of size less than $k$ cannot include a hyperedge, we have
\begin{align}\label{coeffip}
a_i^n=\binom{n}{i}, i=0,1,\cdots,k-1.
\end{align}

In what follows, we first consider the log-concavity of independence polynomials of antiregular $3$-hypergraphs and then that of antiregular $k$-hypergraphs. The core proving strategy is based on the induction on the number of vertices of antiregular hypergraphs.

\begin{lemma}\label{coeff3lem}
The coefficients $a_3^{2n}$ and $a_4^{2n}$ of the independence polynomial $I(\bar{A}_{2n};x)$ of the antiregular $3$-hypergraph $\bar{A}_{2n}$ are given by
\begin{align}
a_3^{2n}=\frac{1}{6}n(n-1)(4n+1),~~ a_4^{2n}=\frac{1}{6}n^2(n-1)(n-2).
\end{align} 
\end{lemma}

\begin{proof}
By \eqref{k3ipeq1}, we get
\begin{align*}
&a_3^{2n}=3\binom{n+1}{3}+\binom{n}{3}=\frac{1}{6}n(n-1)(4n+1),\\
&a_4^{2n}=3\binom{n+1}{4}+\binom{n}{4}=\frac{1}{6}n^2(n-1)(n-2).
\end{align*}
\end{proof}

\begin{theorem}\label{logc3thm}
The independence polynomials of antiregular $3$-hypergraphs $A_n$ and $\bar{A}_n$ are log-concave.
\end{theorem}

\begin{proof}
It is easy to see that $I(A_1;x)=1+x$ and $I(A_2;x)=(1+x)^2$ are log-concave. Suppose $I(A_i;x)$ and $I(\bar{A}_i;x)$ are log-concave for $i=1,2,\cdots,m$. Now we show that $I(A_{m+1};x)$ and $I(\bar{A}_{m+1};x)$ are log-concave.

If $m+1=2n$, by Theorem \ref{recthm}, $I(\bar{A}_{2n};x)=(1+x)I(A_{2n-1};x)$ is log-concave since the product of two log-concave polynomials is log-concave.

If $m+1=2n+1$, by Theorem \ref{recthm}, $I(A_{2n+1};x)=I(\bar{A}_{2n};x)+x(2nx+1)$, which means
\begin{align*}
a_i^{2n+1}=a_i^{2n},~~i\geq3.
\end{align*}
Since $I(\bar{A}_{2n};x)$ is log-concave, then for $i\geq 4$,
\begin{align*}
(a_i^{2n+1})^2-a_{i-1}^{2n+1}a_{i+1}^{2n+1}=(a_i^{2n})^2-a_{i-1}^{2n}a_{i+1}^{2n}\geq 0.
\end{align*}
Note that $(a_1^{2n+1})^2-a_{0}^{2n+1}a_{2}^{2n+1}>0$ because of \eqref{coeffip}, therefore, we only have to check the nonnegativeness of $(a_i^{2n+1})^2-a_{i-1}^{2n+1}a_{i+1}^{2n+1}$ for $i=2,3$.

By Lemma \ref{coeff3lem} and \eqref{coeffip},
\begin{align*}
(a_2^{2n+1})^2-a_{1}^{2n+1}a_{3}^{2n+1}&=\left(\binom{2n+1}{2}\right)^2-\binom{2n+1}{1}\frac{1}{6}n(n-1)(4n+1)\\
&=\frac{1}{6}n(n+1)(2n+1)(8n+1)>0,\\
(a_3^{2n+1})^2-a_{2}^{2n+1}a_{4}^{2n+1}&=\left(\frac{1}{6}n(n-1)(4n+1)\right)^2-\binom{2n+1}{2}\frac{1}{6}n^2(n-1)(n-2)\\
&=\frac{1}{36}n^2(n-1)(n+1)(4n^2+6n-1)\geq 0
\end{align*}
Therefore, we complete this proof by induction.
\end{proof}

\begin{lemma}\label{coeffklem}
The coefficients $a_k^{2n}$ and $a_{k+1}^{2n}$ of the independence polynomial $I(\bar{A}_{2n};x)$ of the antiregular $k$-hypergraph $\bar{A}_{2n}$ are given by
\begin{align}
a_k^{2n}=\sum_{i=\lfloor \frac{k+1}{2}\rfloor}^n \binom{2i-1}{k-1},~~ a_{k+1}^{2n}=\sum_{i=\lfloor \frac{k+1}{2}\rfloor}^{n-1} \binom{2i-1}{k-1}(n-i).
\end{align} 
\end{lemma}

\begin{proof}
By Theorem \ref{recthm}, we have
\begin{equation}\label{lemkeq1}
\begin{aligned}
I(\bar{A}_{2n};x)&=(1+x)I(A_{2n-1};x)\\
&=(1+x)\left[I(\bar{A}_{2n-2};x)+\sum_{i=1}^{k-1}\binom{2n-2}{i-1}x^i\right]\\
&=(1+x)I(\bar{A}_{2n-2};x)+\binom{2n-2}{k-2}x^k+\sum_{i=1}^{k-1}\binom{2n-1}{i-1}x^i.
\end{aligned}
\end{equation}
From \eqref{lemkeq1}, we get the recurrence relations
\begin{align}\label{lemkeq2}
a_k^{2n}=a_k^{2n-2}+a_{k-1}^{2n-2}+\binom{2n-2}{k-2}=a_k^{2n-2}+\binom{2n-1}{k-1},
\end{align}
where we have used \eqref{coeffip} and $\binom{2n-2}{k-1}+\binom{2n-2}{k-2}=\binom{2n-1}{k-1}$. By induction on \eqref{lemkeq2}, we find
\begin{align*}
a_k^{2n}&=\begin{cases}
a_k^{k+1}+\sum_{i=\frac{k+3}{2}}^n\binom{2i-1}{k-1},~~~k~\hbox{is odd}\\
a_k^k+\sum_{i=\frac{k+2}{2}}^n\binom{2i-1}{k-1},~~~k~\hbox{is even}
\end{cases}\\
&=\begin{cases}
\sum_{i=\frac{k+1}{2}}^n\binom{2i-1}{k-1},~~~k~\hbox{is odd}\\
\sum_{i=\frac{k}{2}}^n\binom{2i-1}{k-1},~~~k~\hbox{is even}
\end{cases}\\
&=\sum_{i=\lfloor \frac{k+1}{2}\rfloor}^n \binom{2i-1}{k-1}.
\end{align*}
Similarly, from \eqref{lemkeq1}, we obtain
\begin{align*}
a_{k+1}^{2n}&=a_{k+1}^{2n-2}+\sum_{i=\lfloor \frac{k+1}{2}\rfloor}^{n-1} \binom{2i-1}{k-1}\\
&=\begin{cases}
a_{k+1}^{k+3}+\sum_{i=\frac{k+3}{2}}^{n-1}\sum_{i=\frac{k+1}{2}}^j\binom{2i-1}{k-1},~~~k~\hbox{is odd}\\
a_{k+1}^{k+2}+\sum_{i=\frac{k+2}{2}}^{n-1}\sum_{i=\frac{k}{2}}^j\binom{2i-1}{k-1},~~~k~\hbox{is even}
\end{cases}\\
&=\begin{cases}
\sum_{i=\frac{k+1}{2}}^{n-1}\sum_{i=\frac{k+1}{2}}^j\binom{2i-1}{k-1},~~~k~\hbox{is odd}\\
\sum_{i=\frac{k}{2}}^{n-1}\sum_{i=\frac{k}{2}}^j\binom{2i-1}{k-1},~~~k~\hbox{is even}
\end{cases}\\
&=\sum_{i=\lfloor \frac{k+1}{2}\rfloor}^{n-1} \binom{2i-1}{k-1}(n-i).
\end{align*}
\end{proof}

\begin{theorem}\label{logckthm}
The independence polynomials of antiregular $k$-hypergraphs $A_n$ and $\bar{A}_n$ are log-concave.
\end{theorem}

\begin{proof}
It is easy to see that $I(A_i;x)=(1+x)^i$ are log-concave for $1\leq i\leq k-1$. Suppose $I(A_i;x)$ and $I(\bar{A}_i;x)$ are log-concave for $i=1,2,\cdots,m$ $(m\geq k-1)$. Now we show that $I(A_{m+1};x)$ and $I(\bar{A}_{m+1};x)$ are log-concave.

If $m+1=2n$, by Theorem \ref{recthm}, $I(\bar{A}_{2n};x)=(1+x)I(A_{2n-1};x)$ is log-concave since the product of two log-concave polynomials is log-concave.

If $m+1=2n+1$, by Theorem \ref{recthm}, 
\begin{align*}
I(A_{2n+1};x)=I(\bar{A}_{2n};x)+\sum_{i=1}^{k-1}\binom{2n}{i-1}x^i, 
\end{align*}
which means
\begin{align*}
a_i^{2n+1}=a_i^{2n},~~i\geq k.
\end{align*}
Denote $F_i^{2n+1}:=(a_i^{2n+1})^2-a_{i-1}^{2n+1}a_{i+1}^{2n+1}$. Since $I(\bar{A}_{2n};x)$ is log-concave, then for $i\geq k+1$,
\begin{align*}
F_i^{2n+1}=(a_i^{2n})^2-a_{i-1}^{2n}a_{i+1}^{2n}\geq 0.
\end{align*}
Note that $F_i^{2n+1}\geq 0$ for $1\leq i \leq k-2$ because of \eqref{coeffip}, therefore, we only have to check the nonnegativeness of $F_i^{2n+1}$ for $i=k-1,k$.

By Lemma \ref{coeffklem} and \eqref{coeffip},
\begin{align*}
F_{k-1}^{2n+1}&=\left(\binom{2n+1}{k-1}\right)^2-\binom{2n+1}{k-2}a_k^{2n}\\
&>\left(\binom{2n+1}{k-1}\right)^2-\binom{2n+1}{k-2}\binom{2n+1}{k}\geq 0.
\end{align*}
\begin{align*}
F_{k}^{2n+1}&=\left(\sum_{i=\lfloor \frac{k+1}{2}\rfloor}^n \binom{2i-1}{k-1}\right)^2-\binom{2n+1}{k-1}\left(\sum_{i=\lfloor \frac{k+1}{2}\rfloor}^{n-1} \binom{2i-1}{k-1}(n-i)\right)\\
&\geq \left(\sum_{i=\lfloor \frac{k+1}{2}\rfloor}^n \binom{2i-1}{k-1}\right)^2-\binom{2n+1}{k-1}\left(\sum_{i=\lfloor \frac{k+1}{2}\rfloor}^{n-1} \binom{2i-1}{k-1}\sum_{i=\lfloor \frac{k+1}{2}\rfloor}^{n-1}(n-i)\right)\\
&=\left(\sum_{i=\lfloor \frac{k+1}{2}\rfloor}^{n} \binom{2i-1}{k-1}\right)^2-\binom{n+1-\lfloor \frac{k+1}{2}\rfloor}{2}\binom{2n+1}{k-1}\sum_{i=\lfloor \frac{k+1}{2}\rfloor}^{n-1} \binom{2i-1}{k-1}\\
&=\left(\sum_{i=\lfloor \frac{k+1}{2}\rfloor}^{n-1} \binom{2i-1}{k-1}\right)^2+\left(2\binom{2n-1}{k-1}-\binom{n+1-\lfloor \frac{k+1}{2}\rfloor}{2}\binom{2n+1}{k-1}\right)\times \\
&\sum_{i=\lfloor \frac{k+1}{2}\rfloor}^{n-1} \binom{2i-1}{k-1}+\left(\binom{2n-1}{k-1}\right)^2=:\tilde{F}.
\end{align*}
We observe that $\tilde{F}$ is a quadratic form of $\sum_{i=\lfloor \frac{k+1}{2}\rfloor}^{n-1} \binom{2i-1}{k-1}$ and the quadratic form
\begin{align*}
r^2+\left(2\binom{2n-1}{k-1}-\binom{n+1-\lfloor \frac{k+1}{2}\rfloor}{2}\binom{2n+1}{k-1}\right)r+\left(\binom{2n-1}{k-1}\right)^2
\end{align*} 
has two positive roots,
\begin{align*}
r_{\pm}=\frac{1}{2}\left(m-2\pm\sqrt{m^2-4m}\right)\binom{2n-1}{k-1},
\end{align*}
where $m=\frac{2n(2n+1)}{(2n-k+1)(2n-k+2)}\binom{n+1-\lfloor \frac{k+1}{2}\rfloor}{2}$. By induction hypothesis, $F_k^{2n-1}\geq 0$ gives
\begin{align*}
\sum_{i=\lfloor \frac{k+1}{2}\rfloor}^{n-1} \binom{2i-1}{k-1}\geq \binom{2n-1}{k-1}\left(\sum_{i=\lfloor \frac{k+1}{2}\rfloor}^{n-2} \binom{2i-1}{k-1}(n-1-i)\right).
\end{align*}
For $n\geq \lfloor \frac{k+1}{2}\rfloor+3$, 
\begin{align*}
\sum_{i=\lfloor \frac{k+1}{2}\rfloor}^{n-2} \binom{2i-1}{k-1}(n-1-i)\geq m,
\end{align*}
which is because the left side is at least $O(n^3)$ while the right side is $O(n^2)$. Thus, $\tilde{F}\geq 0$, which implies $F_k^{2n+1}\geq 0$. For $n=\lfloor \frac{k+1}{2}\rfloor+1$ or $n=\lfloor \frac{k+1}{2}\rfloor+2$, by straightforward calculations, we get
\begin{align*}
F_k^{2\lfloor \frac{k+1}{2}\rfloor+3}&=\begin{cases}
\frac{k^2}{360}\left(7k^4+30k^3+145k^2+330k+568\right),~~&k~\hbox{is odd}\\
\frac{1}{24}\left(5k^4+6k^3+19k^2+18k+24\right),&k~\hbox{is even}
\end{cases}\\
&>0,\\
F_k^{2\lfloor \frac{k+1}{2}\rfloor+5}&=\begin{cases}
\frac{11}{302400}k^{10}+\ldots,~~&k~\hbox{is odd}\\
\frac{1}{960}k^8+\ldots,&k~\hbox{is even}
\end{cases}\\
&>0.
\end{align*}
Therefore, we complete this proof by induction.
\end{proof}

\section{Threshold of antiregular $k$-hypergraphs}\label{sec4}

Chv{\'a}tal and Hammer \cite{chvtal1977aggregation} introduced {\em threshold graphs} as the graphs with the following property: a simple graph $G=(V,E)$ is called {\em threshold} if there exist a labeling $c$ of $V$ and a threshold $\tau$ such that $X\subseteq V$ is {\em stable}\footnote{It is now commonly called {\em independent}.}, i.e., any two vertices in $X$ are not adjacent, if and only if $\sum_{x\in X}c(x)\leq \tau$. Note that there are several equivalent definitions of threshold graphs, readers may refer to \cite{golumbic1976threshold,golumbic2004algorithmic,henderson1977graph}. Golumbic \cite{golumbic2004algorithmic} suggested to generalise the notion of threshold graphs to that of threshold hypergraphs and to study their properties. The aim of this section is to show that antiregular $k$-hypergraphs and a broader class of hypergraphs, $\{0,1\}$-constructable $k$-hypergraphs, are $T2$-threshold.

Recall the definitions of threshold hypergraphs proposed by Golumbic \cite{golumbic2004algorithmic}.

\begin{definition}
Let $H=(V,E)$ be a $k$-hypergraph, we consider the following properties:

(T1) There exist a labeling $c$ of $V$ and a threshold $\tau$ such that, for any vertex subset $X\subseteq V$, $X$ contains a hyperedge if and only if $\sum_{x\in X}c(x)>\tau$.

(T2) There exist a labeling $c^\prime$ of $V$ and a threshold $\tau^\prime$ such that, for any vertex subset $X^\prime\subseteq V$ of size $k$, $X^\prime \in E$ if and only if $\sum_{x\in X^\prime}c^\prime(x)>\tau^\prime$.

(T3) For $x,y\in V$, define $x\ll y$ if $x$ can be replaced by $y$ in any hyperedge, i.e., if for any $\{x_1,x_2,\ldots,x_{k-1}\}\in \vert V\setminus \{x,y\}\vert^{k-1}$, $\{x,x_1,\ldots,x_{k-1}\}\in E$ $\implies$ $\{y,x_1,\ldots,x_{k-1}\}\in E$. Then, for any $x,y\in V$, either $x\ll y$ or $y\ll x$ or both holds.

The hypergraph $H$ is called {\em $Ti$-threshold} if it satisfies (Ti), where $i=1,2,3$.
\end{definition}

It is easy to see that $(T1)\implies (T2)\implies (T3)$ and the reversed implications are true for $k=2$. However, for $k\geq 3$, neither $(T3)\implies (T2)$ nor $(T2)\implies (T1)$ holds. Some counterexample are outlined in \cite{golumbic2004algorithmic,reiterman1985threshold}.

\begin{definition}
A $k$-hypergraph is called {\em $\{0,1\}$-constructable} if it can be constructed by the operations of adding an isolated vertex (connected to no other vertices) and adding a dominating vertex (connected to every other $k-1$ vertices) in some order.
\end{definition}

\begin{figure}[!h]
\includegraphics[width=4cm]{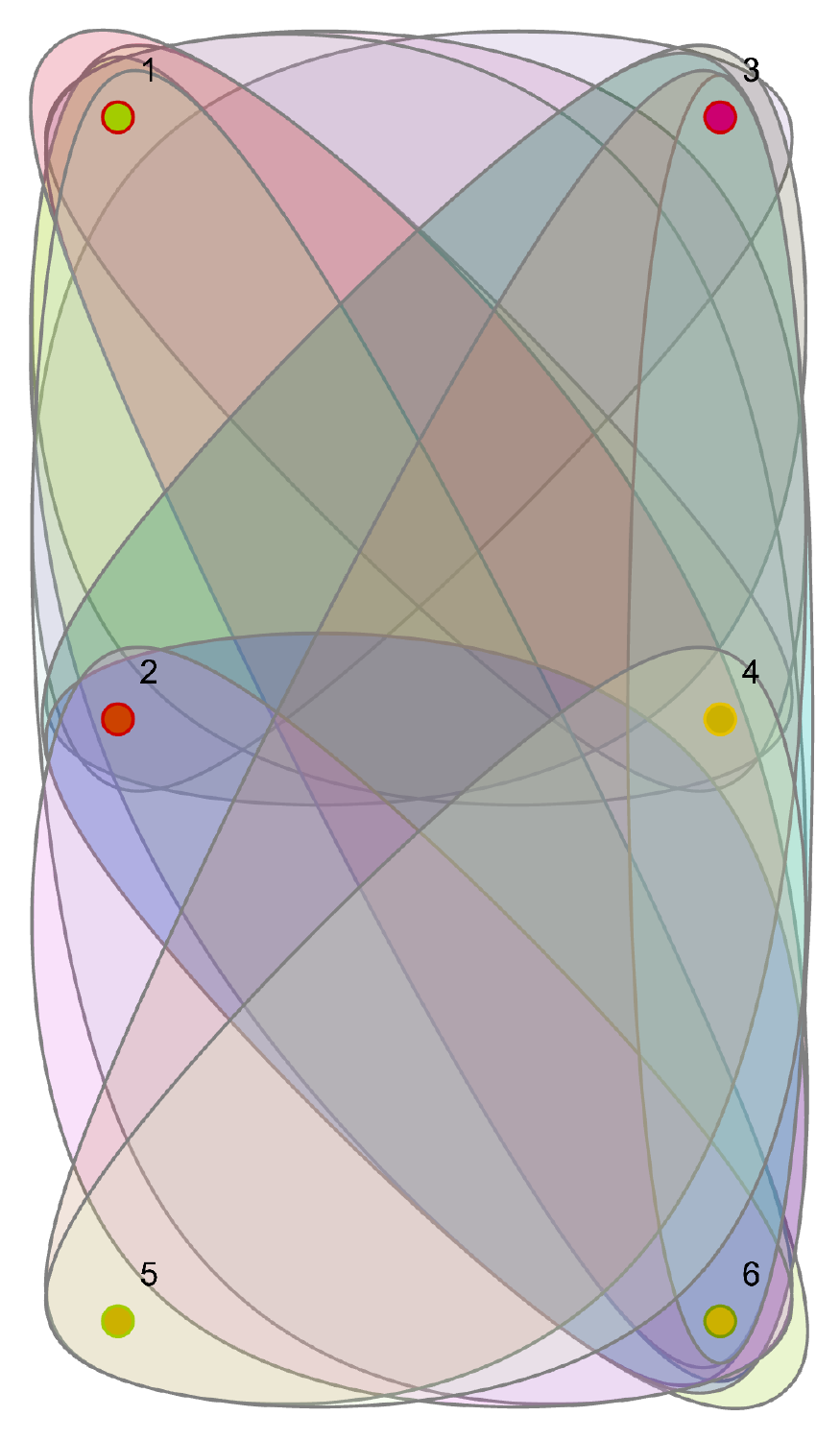}
\hspace{6mm}
\includegraphics[width=4cm]{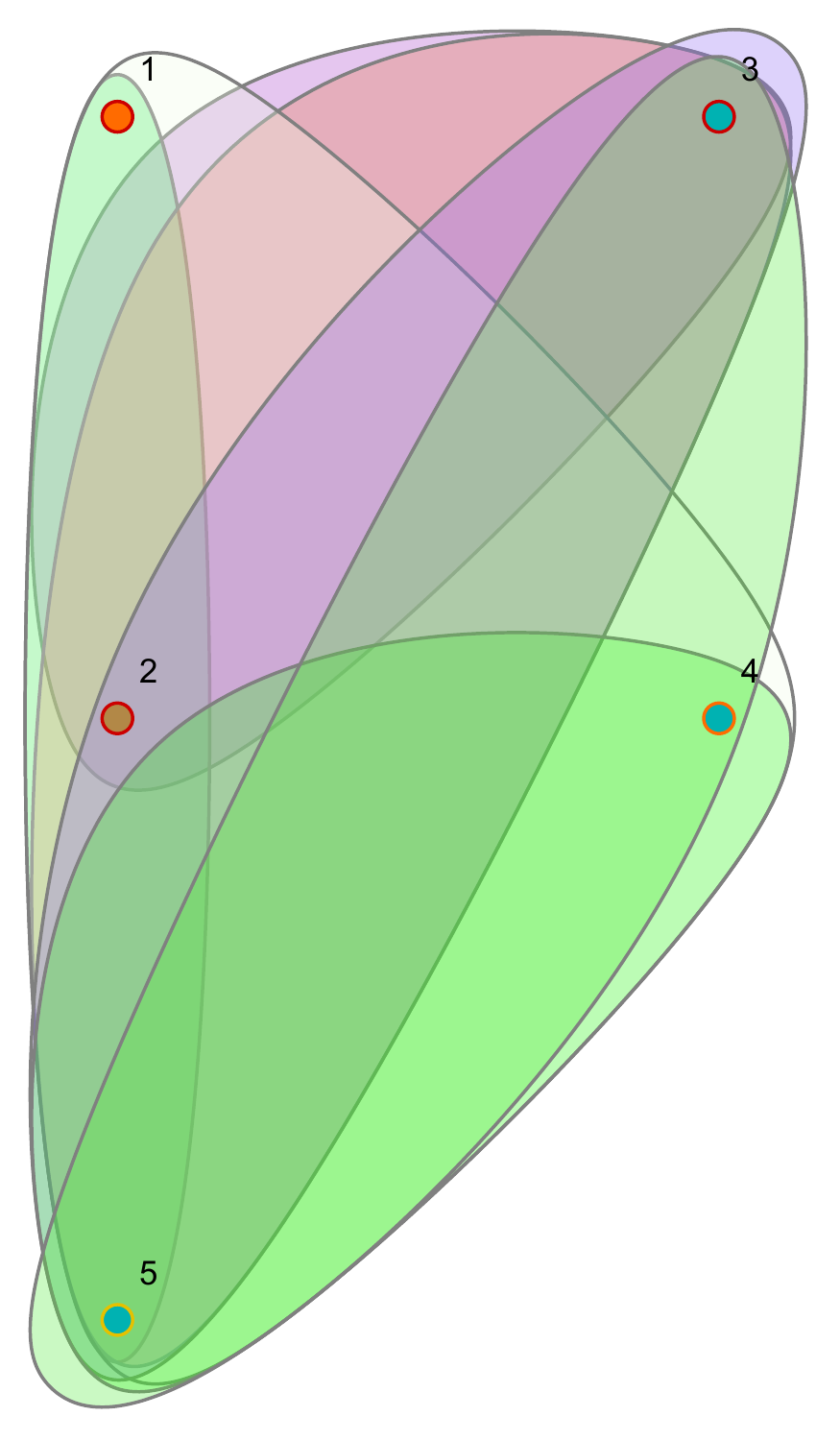}
\caption{Two $\{0,1\}$-constructable $3$-hypergraphs corresponding to the binary building strings $001101$ (left) and $00101$ (right), respectively.}
\label{fig1}
\end{figure}

Not all $k$-hypergraphs are $\{0,1\}$-constructable. For example, the $3$-hypergraph $H=([6],\{\{1,2,3\},\{3,4,5\},\{1,5,6\}\})$ is not $\{0,1\}$-constructable, where $[6]=\{1,2,\ldots,6\}$. Obviously, all antiregular $k$-hypergraphs are $\{0,1\}$-constructable. Figure \ref{fig1} shows two examples of $\{0,1\}$-constructable $3$-hypergraphs, where the binary building strings are $001101$ and $00101$ (antiregular).

For simplicity, we denote by $b^n$ the binary building string of a $\{0,1\}$-constructable $k$-hypergraph with $n$ vertices, $c$ the labeling of $b^n$, $b_i^n$ the $i$-th bit of $b^n$ ($1\leq i\leq n$), $b^{n,j}:=\{i\in [n]\vert b_i^n=j\}$ the set of all isolated vertices ($j=0$) or dominating vertices ($j=1$), $b^{n,j}_{m}$ the last $m$ vertices in $b^{n,j}$, i.e.,
\begin{align*}
b^{n,j}_m:=&\{i_1,i_2,\ldots,i_m\in b^{n,j}\vert~  i_1<i_2<\cdots<i_m=\max b^{n,j},\\
& \forall s\in (i_1,i_m)\setminus\{i_2,\ldots,i_{m-1}\}, b^{n}_s=1-j\}, j=0,1,
\end{align*} 
$b^{n,0}_{c,m}$ the subset of $b^{n,0}$ that $c$ takes smallest $m$ values on it, i.e.,
\begin{align*}
b^{n,0}_{c,m}:=&\{i_1,i_2,\ldots,i_m\in b^{n,0}\vert~  c(i_1)\leq c(i_2)\leq\cdots\leq c(i_m),\\
& \forall s\in b^{n,0}\setminus\{i_1,\ldots,i_{m}\}, c(s)\geq c(i_m)\}.
\end{align*}
For example, let $b^9=001010001$ with $c=[32,32,48,24,56,4,6,7,102]$, then $b^{7,0}=\{1,2,4,6,7,8\}$, $b^{7,1}=\{3,5,9\}$, $b^{7,0}_3=\{6,7,8\}$, $b^{7,1}_2=\{5,9\}$, $b^{7,0}_{c,2}=\{6,7\}$.

Now we outline the algorithm of defining the labeling and threshold for $\{0,1\}$-constructable $k$-hypergraphs as follows.

\begin{algorithm}
\caption{Construction of labels and thresholds for $\{0,1\}$-constructable $k$-hypergraphs}\label{alg:threshold}
\begin{algorithmic}
\Require Binary building string $b^n$ of a $\{0,1\}$-constructable $k$-hypergraph $n$ vertices $(n\geq k\geq 2)$
\Ensure A labeling $c:[n]\rightarrow \mathbb{Z}$ and a threshold $\tau\in\mathbb{Z}$
\State Suppose the first $s$ bits of $b^n$ are $0$'s and $b^n_{s+1}=1$, where $s\geq k-1$.
\begin{itemize} 
\item If $b^n=0[s]1$, we define a labeling $c:[n]\rightarrow \mathbb{Z}$ as
\begin{align}\label{algeq1}
c(i)=\begin{cases}
2,~~ &\hbox{if}~~1\leq i<n\\ 
3,~~ &\hbox{if}~~i=n
\end{cases} 
\end{align}
and set a threshold $\tau=2k$.
\end{itemize}

\State Let $b^m$ be the restriction of $b^n$ on the first $m$ bits, $c^\prime:[m]\rightarrow \mathbb{Z}$ and $\tau^\prime$ its labeling and threshold. Set $m=s+1$ and continue the following two steps until $m+1=n$.
\begin{itemize}
\item If $b^n_{m+1}=1$, then we define a new labeling $c:[m+1]\rightarrow \mathbb{Z}$ as
\begin{align}\label{algeq2}
c(i)=\begin{cases}
c^\prime(i),~~&\hbox{if}~~1\leq i\leq m\\
\tau^\prime+1-\sum_{j\in b^{m,0}_{c^\prime, k-1}}c^\prime(j), &\hbox{if}~~i=m+1
\end{cases}
\end{align}
and set a new threshold $\tau=\tau^\prime$.

\item If $b^n_{m+1}=0$, then we define a new labeling $c:[m+1]\rightarrow \mathbb{Z}$ as
\begin{align}\label{algeq3}
c(i)=\begin{cases}
2c^\prime(i),~~&\hbox{if}~~1\leq i\leq m\\
2\tau^\prime+1-\sum_{j\in b^{m,1}_{k-1}}2c^\prime(j), &\hbox{if}~~i=m+1~~\hbox{and}~~\vert b^{m,1}\vert\geq k-1\\
2\tau^\prime+1-\sum_{j\in b^{m,1}}2c^\prime(j)-\sum_{j=1}^{k-1-\vert b^{m,1}\vert}2c^\prime(j), &\hbox{if}~~i=m+1~~\hbox{and}~~\vert b^{m,1}\vert< k-1
\end{cases}
\end{align}
and set a new threshold $\tau=2\tau^\prime+1$.

\end{itemize}

\end{algorithmic}
\end{algorithm}

Next we will show that the labeling and threshold defined in Algorithm \ref{alg:threshold} for a $\{0,1\}$-constructable $k$-hypergraph can make it $T2$ thresholdable.

\begin{theorem}\label{thresholdthm}
All $\{0,1\}$-constructable $k$-hypergraphs are $T2$-threshold, which are $T3$-threshold as well. Therefore, all antiregular $k$-hypergraphs are $T2$- and $T3$-threshold.
\end{theorem}

\begin{proof}
By \eqref{algeq1} in Algorithm \ref{alg:threshold}, it is easy to check that $b^n=0[s]1$ is $T2$-threshold.

Suppose $b^n$ is $T2$-threshold for $n\leq m$. We assume the labeling and threshold for $b^m$ are $c^\prime$ and $\tau^\prime$, respectively, and construct the new labeling $c$ and threshold $\tau$ for $b^{m+1}$ by Algorithm \ref{alg:threshold}. Now we show $b^{m+1}$ is $T2$-threshold with respect to $c$ and $\tau$.

Let $S=\{i_1,i_2,\ldots,i_k\}$ $(i_1<i_2<\ldots<i_k)$ be a vertex subset of size $k$.

Case I: $b^{m+1}_{m+1}=1$

(a) If $i_k<m+1$, then by induction hypothesis and \eqref{algeq2},
\begin{align*}
\sum_{j=1}^k c(i_j)=\sum_{j=1}^k c^\prime(i_j)>\tau^\prime=\tau\iff S~\hbox{is a hyperedge.}
\end{align*} 

(b) If $i_k=m+1$, then $S$ is a hyperedge. By Lemma \ref{dec0lem},
\begin{align*}
\sum_{j=1}^k c(i_j)\geq c(m+1)+\sum_{j\in b^{m,0}_{c^\prime, k-1}}c^\prime(j)=\tau^\prime+1>\tau^\prime=\tau.
\end{align*}

Case II: $b^{m+1}_{m+1}=0$

(a) If $i_k<m+1$, then by induction hypothesis and \eqref{algeq3},
\begin{align*}
\sum_{j=1}^k c(i_j)=2\sum_{j=1}^k c^\prime(i_j)\geq 2(\tau^\prime+1)=\tau+1>\tau\iff S~\hbox{is a hyperedge.}
\end{align*} 

(b) If $i_k=m+1$, then $S$ is a non-hyperedge. By Lemma \ref{dec0lem},
\begin{align*}
\sum_{j=1}^k c(i_j)&\leq 
\begin{cases}
c(m+1)+2\sum_{j\in b^{m,1}_{k-1}}c^\prime(j),&\hbox{if}~~\vert b^{m,1}\vert\geq k-1\\
c(m+1)+2\sum_{j\in b^{m,1}}c^\prime(j)+2\sum_{j=1}^{k-1-\vert b^{m,1}\vert}c^\prime(j), &\hbox{if}~~\vert b^{m,1}\vert< k-1
\end{cases}\\
&=2\tau^\prime+1=\tau.
\end{align*}
Thus, we complete this proof.
\end{proof}

Let $b^n$ be the binary building string of a $\{0,1\}$-constructable hypergraph, an interval $[i_1,i_2]$ $(1\leq i_1\leq i_2\leq n)$ is called a {\em $i$-interval} if $b^n_j=i$ for $i_1\leq j\leq i_2$ and $b^n_{i_1-1}=b^n_{i_2+1}=1-i$ (if $i_1-1\geq 1$ and/or $i_2+1\leq n$), where $i=0,1$. Further, $[i_1,i_2]$ is called {\em trivial} if $i_1=i_2$ and {\em non-trivial} otherwise. For instance, $b^n=00110001011$, its $0$-intervals are $[1,2]$, $[5,7]$, $[9,9]$ (trivial), and $1$-intervals are $[3,4]$, $[8,8]$ (trivial), $[10,11]$.

\begin{lemma}\label{dec0lem}
Let $c$ and $\tau$ be the labeling and threshold of $b^n$ as constructed in Algorithm \ref{alg:threshold}, where $[1,s]$ $(s\geq k-1)$ is a $0$-interval of $b^n$.

(1) Let $[i_1,i_2]$ and $[i_3,i_4]$ $(i_2<i_3)$ be any two $1$-intervals of $b^n$, $[i_5,i_6]$ any non-trivial $1$-interval of $b^n$, then
\begin{align}
&c(i)<c(j),~~\hbox{if}~~i\in[i_1,i_2], j\in [i_3,i_4]\label{inc1eq1}\\
&c(i)=c(j),~~\hbox{if}~~i,j\in [i_5,i_6], i<j\label{inc1eq2}
\end{align}
That is, from the left to right, the labels are increasing for vertices in different $1$-intervals and keep the same for vertices in the same $1$-intervals.

(2) Let $[i_1,i_2]$ and $[i_3,i_4]$ $(i_1>s, i_2<i_3)$ be any two $0$-intervals of $b^n$, $[i_5,i_6]$ any non-trivial $0$-interval of $b^n$, then
\begin{align}
&c(i)=c(j),~~\hbox{if}~~i,j\in [1,s]\label{dec0eq1}\\
&c(i)>c(j),~~\hbox{if}~~i\in[i_1,i_2], j\in [i_3,i_4]\label{dec0eq2}\\
&c(i)<c(j),~~\hbox{if}~~i,j\in [i_5,i_6], i<j\label{dec0eq3}
\end{align}
That is, from the left to right, the labels are decreasing for vertices in different $0$-intervals and increasing for vertices in the same $0$-intervals except $[1,s]$.

Therefore, labels of dominating vertices are always larger than those of isolated vertices.
\end{lemma}

\begin{proof}
(1) Let $b^{m+2}=0[s]1\ldots 011$ be the restriction of $b^n$, $c^{l},\tau^{l}$ the labeling and threshold generated by Algorithm \ref{alg:threshold} for $b^{m+l}$, $l=0,1,2$. Then
\begin{align*}
c^2(m+2)&=\tau^1+1-\sum_{j\in b^{m+1,0}_{c^1,k-1}}c^1(j)\\
&=\tau^0+1-\sum_{j\in b^{m,0}_{c^0,k-1}}c^0(j)\\
&=c^1(m+1)=c^2(m+1).
\end{align*}
Therefore, \eqref{inc1eq2} holds by induction.

Let $b^{m+3}=0[s]1\ldots 101$ be the restriction of $b^n$, $c^{l},\tau^{l}$ the labeling and threshold generated by Algorithm \ref{alg:threshold} for $b^{m+l}$, $l=0,1,2,3$. Then
\begin{align*}
c^3(m+3)&=\tau^2+1-\sum_{j\in b^{m+2,0}_{c^2,k-1}}c^2(j)\\
&>2\tau^1+2-\sum_{j\in b^{m,0}_{c^1,k-1}}2c^1(j)\\
&=2\tau^0+2-\sum_{j\in b^{m,0}_{c^1,k-1}}2c^0(j)\\
&=2c^0(m+1)=c^3(m+1).
\end{align*}
Therefore, \eqref{inc1eq1} holds by induction.

(2) \eqref{dec0eq1} follows from \eqref{algeq1} and the labels of the first $s$ bits are either kept the same or doubled simultaneously in all iterations of Algorithm \ref{alg:threshold}.

Let $b^{m+2}=0[s]1\ldots 100$ be the restriction of $b^n$, $c^{l},\tau^{l}$ the labeling and threshold generated by Algorithm \ref{alg:threshold} for $b^{m+l}$, $l=0,1,2$. Then
\begin{align*}
c^2(m+2)&=\begin{cases}
2\tau^1+1-\sum_{j\in b^{m+1,1}_{k-1}}2c^1(j), &\hbox{if}~\vert b^{m+1,1}\vert\geq k-1\\
2\tau^1+1-\sum_{j\in b^{m+1,1}}2c^1(j)-\sum_{j=1}^{k-1-\vert b^{m+1,1}\vert}2c^1(j), &\hbox{if}~\vert b^{m+1,1}\vert< k-1
\end{cases}\\
&=\begin{cases}
4\tau^0+3-\sum_{j\in b^{m,1}_{k-1}}4c^0(j), &\hbox{if}~\vert b^{m,1}\vert\geq k-1\\
4\tau^0+3-\sum_{j\in b^{m,1}}4c^0(j)-\sum_{j=1}^{k-1-\vert b^{m,1}\vert}4c^0(j), &\hbox{if}~\vert b^{m,1}\vert< k-1
\end{cases}\\
&>\begin{cases}
2(2\tau^0+1-\sum_{j\in b^{m,1}_{k-1}}2c^0(j)), &\hbox{if}~\vert b^{m,1}\vert\geq k-1\\
2(2\tau^0+1-\sum_{j\in b^{m,1}}2c^0(j)-\sum_{j=1}^{k-1-\vert b^{m,1}\vert}2c^0(j)), &\hbox{if}~\vert b^{m,1}\vert< k-1
\end{cases}\\
&=2c^1(m+1)=c^2(m+1).
\end{align*}
Therefore, \eqref{dec0eq3} holds by induction.

Let $b^{m+3}=0[s]1\ldots 010$ be the restriction of $b^n$, $c^{l},\tau^{l}$ the labeling and threshold generated by Algorithm \ref{alg:threshold} for $b^{m+l}$, $l=0,1,2,3$. Then
\begin{align*}
c^3(m+3)&=\begin{cases}
2\tau^2+1-\sum_{j\in b^{m+2,1}_{k-1}}2c^2(j), &\hbox{if}~\vert b^{m+2,1}\vert\geq k-1\\
2\tau^2+1-\sum_{j\in b^{m+2,1}}2c^2(j)-\sum_{j=1}^{k-1-\vert b^{m+2,1}\vert}2c^2(j), &\hbox{if}~\vert b^{m+2,1}\vert< k-1
\end{cases}\\
&<\begin{cases}
2\tau^1-\sum_{j\in b^{m+1,1}_{k-1}}2c^1(j), &\hbox{if}~\vert b^{m+1,1}\vert\geq k-1\\
2\tau^1-\sum_{j\in b^{m+1,1}}2c^1(j)-\sum_{j=1}^{k-1-\vert b^{m+1,1}\vert}2c^1(j), &\hbox{if}~\vert b^{m,1}\vert< k-1
\end{cases}\\
&=\begin{cases}
4\tau^0+2-\sum_{j\in b^{m,1}_{k-1}}4c^0(j), &\hbox{if}~\vert b^{m,1}\vert\geq k-1\\
4\tau^0+2-\sum_{j\in b^{m,1}}4c^0(j)-\sum_{j=1}^{k-1-\vert b^{m,1}\vert}4c^0(j), &\hbox{if}~\vert b^{m,1}\vert< k-1
\end{cases}\\
&=\begin{cases}
2(2\tau^0+1-\sum_{j\in b^{m,1}_{k-1}}2c^0(j)), &\hbox{if}~\vert b^{m,1}\vert\geq k-1\\
2(2\tau^0+1-\sum_{j\in b^{m,1}}2c^0(j)-\sum_{j=1}^{k-1-\vert b^{m,1}\vert}2c^0(j)), &\hbox{if}~\vert b^{m,1}\vert< k-1
\end{cases}\\
&=2c^1(m+1)=c^3(m+1).
\end{align*}
Therefore, \eqref{dec0eq2} holds by induction.
\end{proof}

Since there are no non-trivial $1$-intervals in the binary building strings $b^n$ of antiregular $k$-hypergraphs, the labels of dominating vertices are increasing. Except one non-trivial $0$-interval in $b^n$, the labels of isolated vertices are decreasing afterwards. 

\begin{example}
Given a binary building string of a $\{0,1\}$-constructable $3$-hypergraph $b^{13}=0010100011101$, by implementing Algorithm \ref{alg:threshold}, we get the labeling
\begin{align*}
c=[64,64,96,48,112,8,12,14,204,204,204,-185,401]
\end{align*}
and the threshold $\tau=223$.

Given a binary building string of an antiregular $3$-hypergraph $b^{13}=0010101010101$, by implementing Algorithm \ref{alg:threshold}, we get the labeling
\begin{align*}
c=[64,64,96,48,112,8,168,-60,276,-222,506,-559,1005]
\end{align*}
and the threshold $\tau=223$. It is easy to see that these results are consistent with Lemma \ref{dec0lem}.
\end{example}

Threshold graphs are $\{0,1\}$-constructable \cite{chvatal1973set}, however, $T2$-threshold $k$-hypergraphs are not always $\{0,1\}$-constructable. We illustrate this claim by the following example.

\begin{example}
Let $H=(V,E)$ be a $4$-hypergraph, where the vertex set $V=\{w,-2,-1,0,1,2\}$ and the hyperedge set $E=\{\{w\}\cup \{x_1,x_2,x_3\}\vert~\sum_{i=1}^3x_i>0\}$. By defining a labeling $c:V\rightarrow \mathbb{Z}$,
\begin{align*}
c(v)=\begin{cases}
10,~~&v=w,\\
v,&v\in\{-2,-1,0,1,2\}
\end{cases}
\end{align*}
and setting the threshold $\tau=10$, we can easily verify that $H$ is $T2$-threshold.

Suppose $H$ is $\{0,1\}$-constructable. Let $b^6$ be its binary building string, then its first three bits must be $000$. Since $H$ has $4$ hyperedges, i.e.,
\begin{align*}
E=\{\{w,-2,1,2\},\{w,-1,0,2\},\{w,-1,1,2\},\{w,0,1,2\}\},
\end{align*}
the last three bits of $b^6$ could only be $010$. Therefore, $b^n=000010$, whose vertex-degree sequence is $(3,3,3,3,4,0)$, which contradicts with the vertex-degree sequence $(4,1,2,2,3,4)$ of $H$.
\end{example}

The independence polynomials of antiregular graphs are unique within the family of threshold graphs, however, to the author's best knowledge, it is not known whether the independence polynomials of antiregular $k$-hypergraphs are uniquely determined within the family of $T2$-threshold $k$-hypergraphs. The following example shows two non-isomorphic $k$-hyeprgraphs can have the same independence polynomial.

\begin{example}
Let 
\begin{align*}
&H_1=([5],\{\{1,4,5\},\{2,3,5\},\{2,4,5\},\{3,4,5\}\}),\\
&H_2=([5],\{\{1,2,3\},\{1,3,4\},\{2,3,5\},\{3,4,5\}\})
\end{align*}
be two $3$-hypergraphs (see Figure \ref{fig2}), where $[5]=\{1,2,\ldots,5\}$. Then, their vertex-degree sequences are $[1,2,2,3,4]$ and $[2,2,4,2,2]$, which means $H_1$ and $H_2$ are non-isomorphic. Nonetheless, $H_1$ and $H_2$ share the same independence polynomial,
\begin{align*}
I(H_1;x)=I(H_2;x)=1+5x+10x^2+6x^3+x^4.
\end{align*}

\begin{figure}[!h]
\includegraphics[width=4cm]{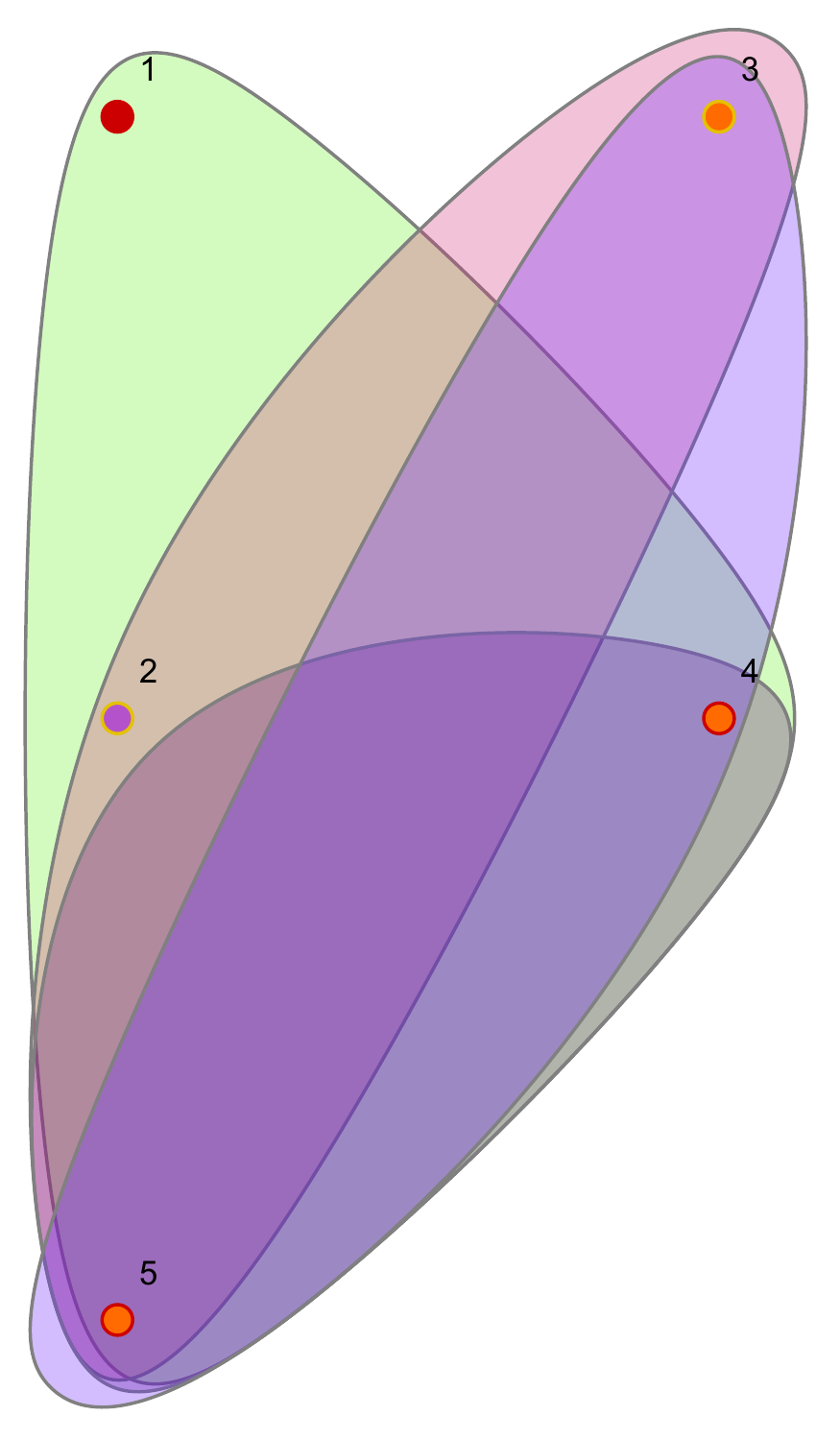}
\hspace{6mm}
\includegraphics[width=4cm]{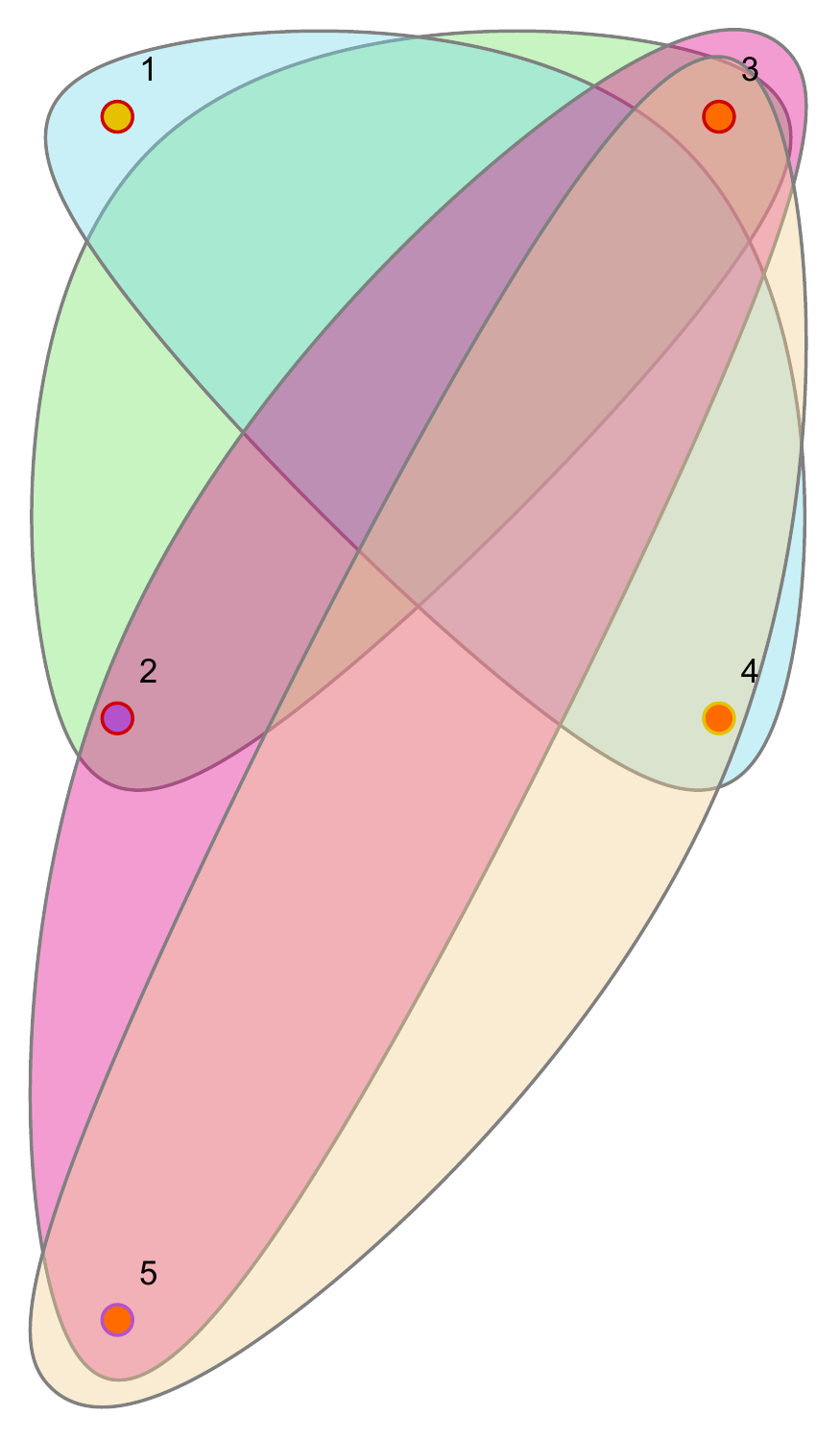}
\caption{Two non-isomorphic $3$-hypergraphs $H_1$ (left) and $H_2$ (right) having the same independence polynomial.}
\label{fig2}
\end{figure}

In this example, $H_1$ is $T2$-threshold, which can be verified by defining a labeling $c_1:[5]\rightarrow \mathbb{Z}$ as $c_1=[-2,-1,0,1,2]$ and setting a threshold $\tau_1=0$. However, $H_2$ is not $T2$-threshold for the following reason. Suppose $H_2$ is $T2$-threshold with respect to a labeling $c_2$ and a threshold $\tau_2$. Since $\{1,3,4\}$ is a hyperedge but not $\{2,3,4\}$, we have
\begin{align*}
c_2(1)+c_2(3)+c_2(4)>\tau_2,~~c_2(2)+c_2(3)+c_2(4)\leq \tau_2,
\end{align*}
which indicates $c_2(1)>c_2(2)$. Since $\{2,3,5\}$ is a hyperedge, we find
\begin{align*}
c_2(1)+c_2(3)+c_2(5)>c_2(2)+c_2(3)+c_2(5)>\tau_2,
\end{align*}
which contradicts with the fact that $\{1,3,5\}$ is not a hyperedge. Therefore, $H_2$ is not $T2$-threshold.
\end{example}

\section{Conclusion}\label{sec5}

Graph properties have been extensively studied in the literature, while hypergraph properties receive much less attention, partly because graphs with some properties cannot be generalised to the hypergraph case naturally. For example, there are many non-equivalent definitions of threshold hypergraphs while they are all equivalent in the graph world by restricting the {\em hyperedge-degree}, i.e., the number of vertices contained in a hyperedge, to $2$. In this paper, we focus on generalising some results of antiregular graphs to those of antiregular $k$-hypergraphs.

In detail, we find the (semi-)closed forms of the independence polynomials of antiregular $k$-hypergraphs. Further, we show that the independence polynomials are log-concave. These results are consistent with Levit and Mandrescu's work \cite{levit2012independence} by setting $k=2$ in the present paper. The third major contribution of this work is that we present an algorithm and prove all $\{0,1\}$-constructable (including antiregular) $k$-hypergraphs are $T2$-threshold. Then, the following implications of properties of $k$-hypergraphs hold,
\begin{align*}
\hbox{antiregular}\implies \{0,1\}\hbox{-constructable}\implies T2\hbox{-threshold}\implies T3\hbox{-threshold}.
\end{align*}
However, the reversed implications are not true for $k\geq 3$. We give an example that $T2$-threshold $\centernot\implies$ $\{0,1\}$-constructable in this paper.

One remaining question here is about the uniqueness of the independence polynomials of antiregular $k$-hypergraphs. Alternatively, a more general question is, {\em given two $T2$-threshold $k$-hypergraphs with the same independence polynomial, are they isomorphic?} These are left for future work.

\begin{acks}[Acknowledgments]
The author thanks Prof. David SUTER for bringing the questions that this paper addresses to my attention and for some helpful discussions of the concepts. The author was partly supported by the Australian Research Council Grant DP200103448. 
\end{acks}

\bibliographystyle{imsart-number} 
\bibliography{references.bib}       

%
%
%

\end{document}